\documentclass[12pt,a4paper]{article}

\usepackage[left=1cm,right=1cm,top=2cm,bottom=2cm]{geometry}
\usepackage[dvips]{graphics}
\usepackage[dvips]{epsfig}
\usepackage{color}
\usepackage{hyperref}
\usepackage{tikz,authblk}

\newtheorem{lemma}{Lemma}
\newtheorem{thm}[lemma]{Theorem}
\newtheorem{cor}[lemma]{Corollary}
\newtheorem{prop}[lemma]{Proposition}
\newtheorem{defn}{Definition}
\newtheorem{rem}{Remark}
\newtheorem{conj}{Conjecture}
\newtheorem{question}{Question}

\newcommand{\dimo}{\noindent \emph{Proof. }}
\newcommand{\qed}{\\ \rightline{$\Box$ \ \ \ \ \ \ \ \ \ \ \ \ \ \ \ }\\}
\newcommand{\e}{\varepsilon}
\newcommand{\G}{\Gamma}
\newcommand{\g}{\gamma}

\usepackage{latexsym}
\usepackage{amsfonts}
\usepackage{amssymb}
\usepackage{amsmath}
\usepackage[english]{babel}

 \title{Compact $4$-manifolds admitting special handle decompositions}

\author[1] {Maria Rita Casali}
\author[2] {Paola Cristofori}

\affil[1] {Department of Physics, Mathematics and Computer Science, University of Modena and Reggio Emilia \ \ \ \ \ \ \ \ \ \ \   Via Campi 213 B, I-41125 Modena (Italy), casali@unimore.it}

\affil[2] {Department of Physics, Mathematics and Computer Science, University of Modena and Reggio Emilia,  \ \ \ \ \ \ \ \ \ \ \   Via Campi 213 B, I-41125 Modena (Italy), paola.cristofori@unimore.it}

\begin{document}
\maketitle

\abstract {In this paper we study colored triangulations of compact PL $4$-manifolds with empty or connected boundary which induce handle decompositions lacking in 1-handles or in 1- and 3-handles, thus 
facing also the problem, posed by Kirby, of the existence of {\it special handlebody decompositions} for any simply-connected closed PL $4$-manifold. In particular, we detect a class of compact simply-connected PL $4$-manifolds 
with empty or connected boundary, which admit such decompositions and, therefore, can be represented by (undotted) framed links. 
Moreover, this class includes any compact simply-connected PL $4$-manifold with empty or connected boundary having colored triangulations that minimize the 
combinatorially defined PL invariant {\em regular genus, gem-complexity} or {\em G-degree} among all such manifolds with the same second Betti number.
}  
 \endabstract
 


\bigskip

\section{Introduction}
\bigskip

The question of whether or not any simply-connected closed PL 4-manifold admits a handle decomposition lacking in 1-handles and 3-handles (i.e. a {\it special handlebody decomposition}, according to \cite[Section~3.3]{[M]}) is long standing: it first appears as problem n. 50 in Kirby's list of problems in low dimensional manifold  topology, edited in 1977  (\cite{Kirby_1977}), it is reported unaltered as Problem 4.18 in the updated list by the same author eighteen years later  (\cite{Kirby_1995}), and so far Kirby's comment about  ``no progress"  continues to be valid. 
% it appears as problem 4.18 in the Kirby's list of problems in low dimensional topology edited in 1995 \cite{Kirby_1995} (and as problem n.50 in the former list by the same author, in 1977: \cite{Kirby_1977}) and 

%Actually, Kirby's problem  presents also an intermediate step with respect to special handlebody decompositions:
Actually, Kirby's problem can be formulated also with an intermediate weaker requirement:

\centerline{\it ``Does every simply-connected closed $4$-manifold have} 
\centerline{\it a handlebody decomposition without 1-handles? } 
\centerline{\it Without 1- and 3-handles?"} 

\medskip

Significant investigations on the matter have been performed during the decades, taking into account  simply-connected PL 4-manifolds with connected boundary, too: recall, for example, the important Trace's contributions  \cite{Trace_1980} and \cite{Trace_1982}. In the boundary case, it is known that each contractible 4-manifold different from $\mathbb D^4$ must have 1- or  3-handles, while a contractible 4-manifold requiring 1-handles in any handle decomposition was detected by Casson  (see \cite[Problem 4.18]{Kirby_1995}).
As regards closed 4-manifolds, the so called {\it Dolgachev surface} $E(1)_{2,3}$, exotic copy of $\mathbb{CP}^2 \#_9 \overline{\mathbb{CP}^2}$, has long been conjectured to be a possible counterexample to the % existence ? 
universality of special handlebody decompositions (\cite{Harer-Kas-Kirby}); however, Akbulut not only  proved $E(1)_{2,3}$ to admit a handle decomposition lacking in 1-handles and 3-handles (\cite{Akbulut_1}), but also produced an infinite family of exotic copies of  $\mathbb{CP}^2 \#_9 \overline{\mathbb{CP}^2}$ with the same property (\cite{Akbulut_2}).  

\bigskip
%The present paper is devoted to face the above Kirby's problem, in both its steps and extended to compact 4-manifolds with connected boundary,  too, by making use of the so called {\it colored triangulations} (i.e. pseudosimplicial triangulations whose 4-simplices have vertices injectively labelled by five colors), together with their visualizing tool consisting of {\it edge-colored graphs}. 
The present paper studies both questions posed by the above Kirby's problem, also in the extended setting of compact 4-manifolds with connected boundary, by making use of the so called {\it colored triangulations} (i.e. pseudosimplicial triangulations whose 4-simplices have vertices injectively labelled by five colors), together with their visualizing tool consisting of {\it edge-colored graphs}. 

More precisely, a large class of colored triangulations (or, equivalently, of associated 5-colored graphs), possibly representing all simply-connected closed PL 4-manifolds, is identified, so that the induced handle decompositions turn out to be special.

From the point of view of colored triangulations, it is well known that any compact PL $4$-manifold with empty or connected boundary admits a {\it contracted triangulation}, i.e. a colored triangulation with exactly five vertices, labelled by the five elements of the color set  $\Delta_4=\{0,1,2,3,4\}$\footnote{Actually, the result holds in general dimension $n$, by simply substituting $n+1$ to $5$, and it is usually stated in terms of the associated graphs, by making use of the notion of {\it crystallization}, which gives the name to the whole combinatorial representation theory via edge-colored graphs: see \cite{Casali-Cristofori-Grasselli}, or Section \ref{preliminaries} of the present paper.}. 
In the closed setting, the main result of the present paper - involving both steps of the quoted Kirby's problem - is the following: 

\begin{thm} \label{thm.special-handle-decompositions_closed}
Let $M^4$ be a closed PL 4-manifold and let $K$ be a contracted triangulation of $M^4$.
%a colored triangulation of $M^4$, containing exactly five vertices. 
\begin{itemize}
\item[(a)] If %a colored triangulation of $M^4$ exists, containing exactly five vertices and 
$K$ contains exactly one edge between two pairs of vertices, involving three vertices, then $M^4$ admits a handle decomposition without 1-handles;  
\item[(b)] if $K$ contains 
% a colored triangulation of $M^4$ exists, containing exactly five vertices and 
exactly one edge between three pairs of vertices, involving all five vertices, then $M^4$ admits a special handlebody decomposition. 
\end{itemize}
\end{thm}

%In the connected boundary case, a contracted triangulation of a compact $4$-manifold $M^4$ is actually a contracted triangulation of the associated {\it singular $4$-manifold} $\widehat{M^4}$ obtained from $M^4$ by capping off the boundary $\partial M^4$ with a cone; 
In the connected boundary case, by a {\it contracted triangulation} of a compact $4$-manifold $M^4$ we mean a contracted triangulation of the associated {\it singular $4$-manifold} $\widehat{M^4}$ obtained from $M^4$ by capping off the boundary $\partial M^4$ with a cone;
moreover, color $4$ is usually assumed to be {\it singular}, i.e. to label the only vertex whose link is a closed connected $3$-manifold different from the $3$-sphere $\mathbb S^3$ (namely, the boundary  $\partial M^4$, which is supposed to be non-spherical).  In this setting, we prove the following result:

\begin{thm} \label{thm.special-handle-decompositions_boundary}
Let $M^4$ be a compact PL 4-manifold with non-empty connected boundary and let $K$ be a contracted triangulation of $M^4$ with $4$ as (unique) singular color.
%a colored triangulation of $M^4$, containing exactly five vertices. 
\begin{itemize}
\item[(a)] If %a colored triangulation of $M^4$ exists, containing exactly five vertices and 
$K$ contains exactly one edge between two pairs of non-singular vertices, involving three vertices, then $M^4$ admits a handle decomposition without 1-handles;  
\item[(b)] if $K$ contains 
% a colored triangulation of $M^4$ exists, containing exactly five vertices and 
exactly one edge between three pairs of vertices, involving all five vertices, and so that only one pair involves color $4$, then $M^4$ admits a special handlebody decomposition. 
\end{itemize}
\end{thm}

\bigskip 

According to a well-known literature (see, for example, \cite{[M]} or \cite{[K]}), a framed link (possibly with dotted circles) can be associated to any handle decomposition of a closed %compact 
PL 4-manifold: roughly speaking, dotted circles represent 1-handles, while undotted components % of the framed link 
give the instructions about the attachment of 2-handles. 
Then, a celebrated result by Montesinos and  Laudenbach-Poenaru (\cite{Montesinos} and \cite{Laudenbach-Poenaru}) ensures that there is a unique way to add the (possible) 3-handles and the 4-handle,
so to obtain the closed $4$-manifold. For this reason, given a framed link with $r\ge 0$ dotted circles $(L,c)$, it is usual to denote by $M^3(L,c)$ the 3-manifold obtained from $\#_r(\mathbb S^2 \times \mathbb S^1)$ by Dehn surgery along the undotted components of $(L,c)$, %to denote 
by $M^4(L,c)$ the compact PL 4-manifold  (with boundary $M^3(L,c)$) obtained from $\mathbb D^4$ by attaching 1- and 2- handles according to the framed link $(L, c)$ and - in case  $M^3(L,c)$ is  
PL-homeomorphic to either $\mathbb S^3$ or a connected sum of $s\ge 0$ copies of the orientable $\mathbb S^2$-bundle over $\mathbb S^1$ - %to denote 
by $\overline{M^4(L,c)}$ the closed  PL 4-manifold obtained from $M^4(L,c)$ by %uniquely 
adding $s$ 3-handles and one 4-handle.  

%In the connected boundary case, under the assumption of simply-connectedness, Trace proved a similar property of unique attachment for 3-handles: 
In the connected boundary case, under the assumption of simply-connectedness, Trace proved a uniqueness property for the attachment of 3-handles:

\begin{prop} \label{prop_Trace_3-handles} \ {\rm \cite[Theorem 1]{Trace_1982}} \ 
Let $N_1^4$, $N_2^4$ be compact simply-connected PL 4-manifolds with non-empty connected boundary, so that $N_i^4$  is obtained  from the compact 4-manifold $W_i^4$ ($\forall i \in \{1,2\}$) by adding $s\ge 0$ 3-handles. 
Then,  $N_1^4$ is PL-homeomorphic to  $N_2^4$ if and only if  $W_1^4$ is PL-homeomorphic to  $W_2^4.$
\end{prop}
  
As Trace himself points out, Proposition  \ref{prop_Trace_3-handles} implies {\it ``that the essential structure of such 4-manifolds is contained in a neighborhood of their 2-skeleton".}  
Hence, framed links (possibly with dotted circles), together with the number of required 3-handles, can
unambiguously identify compact simply-connected PL 4-manifolds with non-empty connected boundary.
 
\bigskip 

As a consequence, if  $M^4(L,c)^{(s)}$ denotes the simply connected compact PL 4-manifold obtained from $M^4(L,c)$  by adding $s$ 3-handles,  then both Theorems  \ref{thm.special-handle-decompositions_closed} and  \ref{thm.special-handle-decompositions_boundary}  may be re-stated as follows:

\begin{prop} \label{prop_link} 
\ 
\begin{itemize} 
\item[(a)] If $M^4$ is a closed  PL 4-manifold  (resp. a compact PL 4-manifold with connected non-empty boundary)  % which admits a contracted triangulation 
satisfying the hypothesis of Theorem  \ref{thm.special-handle-decompositions_closed} (a) (resp.  Theorem  \ref{thm.special-handle-decompositions_boundary} (a)), then a framed link $(L, c)$  (with no dotted circle) and an integer $s \ge 0$ exist, so that:  
$$ M^3(L,c)=\#_s(\mathbb S^2 \times \mathbb S^1) \  \text{ \ and \ }  \ \overline{M^4(L,c)}= M^4$$ % $$\text{ \ with  \ }   \overline M^4(L,c)  \text{ \ obtained from  \ }  M^4(L,c) \text{ \ by adding  $s$ 3-handles}$$
$$ \text{(resp. \ }  M^3(L,c)=\#_s(\mathbb S^2 \times \mathbb S^1)\# \partial M^4  \ \text{ \ and \ }  \ M^4(L,c)^{(s)}= M^4 \text{)}.$$ 
\item[(b)]  If $M^4$ is a closed  PL 4-manifold  (resp. a compact PL 4-manifold with connected non-empty boundary)  %which admits a contracted triangulation 
satisfying the hypothesis of Theorem  \ref{thm.special-handle-decompositions_closed} (b) (resp.  Theorem  \ref{thm.special-handle-decompositions_boundary} (b)), then a framed link $(L, c)$  (with no dotted circle) exists, so that:
$$ M^3(L,c)=\mathbb S^3 \ \text{ \ and \ }  \ \overline{M^4(L,c)}= M^4 $$
$$ \text{(resp. \ } M^3(L,c)= \partial M^4  \  \text{ \ and \ }   \ M^4(L,c)= M^4 \text{)}.$$ 
\end{itemize}
\end{prop}

\bigskip

Within crystallization theory, two interesting classes of graphs representing compact simply-connected PL 4-manifolds (with empty or connected boundary) have been recently introduced and studied: they are 
called {\it weak-simple} and {\it simple crystallizations} and are proved to minimize combinatorially defined PL-invariants, i.e. the {\it regular genus} and, in case of simple crystallizations, also the 
{\it gem-complexity} and  the {\it G-degree} (see \cite{Casali-Cristofori-Gagliardi_Rendiconti_TS_2020} and references therein, or Section \ref{sec-weaksimple}).

Since weak simple (resp. simple) crystallizations - by definition - combinatorially visualize contracted colored triangulations which admit exactly one edge between two vertices labelled by colors consecutive in a chosen permutation of the color set (resp. between any pair of vertices\footnote{Note that, as their names suggest, the class of simple crystallization is a subclass of that of weak simple ones; actually, they both are particular cases of even larger classes of crystallizations, representing not simply-connected compact PL 4-manifolds, and still minimizing combinatorially defined PL invariants, i.e. {\it weak semi-simple} and {\it semi-simple} crystallizations. Precise definitions of all these classes may be found in \cite{Casali-Cristofori-Gagliardi_Rendiconti_TS_2020}, together with its references.}), the represented 4-manifold trivially satisfy the hypothesis of  Theorem  \ref{thm.special-handle-decompositions_closed} (b) (in the closed case) or Theorem  \ref{thm.special-handle-decompositions_boundary} (b) (in the connected boundary case). Hence, the following result is a direct consequence of Theorem \ref{thm.special-handle-decompositions_closed} and  Theorem \ref{thm.special-handle-decompositions_boundary}, together with Proposition \ref{prop_link}: 

\begin{thm} \label{thm.simple&weaksimple}
Let $M^4$ be a compact PL 4-manifold with empty or connected boundary. 
If $M^4$ admits a weak simple crystallization, then $M^4$ admits a special handlebody decomposition (or, equivalently, $M^4$ is represented by a (not dotted) framed link  with $\beta_2(M^4)$ components).
\end{thm}

Note that the above result extends to the connected boundary case a property already proved in \cite{Casali-Cristofori-GagliardiJKTR2015} for closed PL 4-manifolds admitting simple crystallizations; moreover, both in the closed and in the boundary case, our statement takes into account also the larger class of weak simple crystallizations. This latter extension is of particular interest with respect to the quoted Kirby's problem, since it is not difficult to prove the existence of (closed) simply-connected PL 4-manifolds which can't admit simple crystallizations, while it is an open question the possible existence of weak simple crystallizations. 
 
For this reason, we formulate the following  

%\begin{conj} \label{conj_weaksimple} Any closed simply-connected PL 4-manifold admits a weak-simple crystallization.
% \end{conj}

\begin{question}\label{q1} Does any closed simply-connected PL 4-manifold admits a weak-simple crystallization?\end{question}

% In virtue of Theorem \ref{thm.simple&weaksimple}, the above conjecture implies the universality of simple handlebody  decompositions for any 
% simply-connected 4-manifold; moreover, since 4-manifolds admitting weak simple crystallizations are proved to be characterized by the equality between their regular genus and twice their second Betti number (Proposition \ref{reg-gen-Betti-number}), Conjecture \ref{conj_weaksimple} implies also the non-existence of exotic copies of $\mathbb S^4$ (and hence the 4-dimensional Poincar\'e conjecture in smooth category) and of $\mathbb{CP}^2$. 

In virtue of Theorem \ref{thm.simple&weaksimple}, a positive answer to the above question implies the universality of simple handlebody  decompositions for any closed simply-connected 4-manifold; moreover, since 4-manifolds admitting weak simple crystallizations are proved to be characterized by the equality between their regular genus and twice their second Betti number (Proposition \ref{reg-gen-Betti-number}), it implies also the non-existence of exotic copies of $\mathbb S^4$ (and hence the 4-dimensional Poincar\'e conjecture in smooth category) and of $\mathbb{CP}^2$.

%Moreover, via Proposition \ref{reg-gen-Betti-number}, it is not difficult to check the equivalence between  Conjecture  \ref{conj_weaksimple} and a deeply investigated conjecture settled in 1986 within crystallization theory: 

Moreover, via Proposition \ref{reg-gen-Betti-number}, it is not difficult to check the equivalence between  Question \ref{q1} % \ref{conj_weaksimple} 
and a deeply investigated conjecture settled in 1986 within crystallization theory:

\begin{conj} \label{conj_additivity} 
{\rm \cite{Ferri-Gagliardi-Grasselli}} The graph-defined PL invariant regular genus is additive with respect to connected sum of closed simply-connected 4-manifolds. 
\end{conj}

Actually, the achievements of the present paper suggest to formulate the following weaker conjecture, still having all the described implications as regards exotic structures, Poincar\'e Conjecture and special handlebody decompositions:   %%% siamo sicure di TUTTE le implicazioni???

\begin{conj} \label{conj_weaker} Any closed simply-connected PL 4-manifold admits a contracted triangulation satisfying the hypothesis of Theorem \ref{thm.special-handle-decompositions_closed} (b).
\end{conj}

Finally, we point out that the described open questions might be faced also via a computational  approach, in virtue of the totally combinatorial nature of the chosen representing tool: see for example \cite{Casali-Cristofori ElecJComb 2015},  yielding an algorithmic procedure for generating and studying a 4-manifolds catalogue via colored graphs.  
% where an algorithmic procedure for generating and studying a 4-manifolds catalogue via colored graphs is presented.  

\bigskip
\bigskip

\section{Preliminaries} \label {preliminaries} 

In this section we will introduce the basic notions of the representation theory of compact PL manifolds by means of regular edge-colored graphs; further details can 
be found in \cite{Ferri-Gagliardi-Grasselli}. 
All manifolds are supposed to be PL and, when not otherwise stated, connected.

\begin{defn}\label{singular manifold} {\em A {\em singular $n$-manifold} is a closed connected $n$-dimensional polyhedron admitting a simplicial triangulation where 
the links of vertices are closed connected $(n-1)$-manifolds\footnote{Note that, as a consequence, the links of all $h$-simplices of the triangulation with $h > 0$ are $(n-h-1)$-spheres.}. 
Vertices whose links are not  $(n-1)$-spheres are called {\em singular}}.
\end{defn}

\begin{rem} \label{correspondence-sing-boundary} {\em If $N$ is a singular $n$-manifold, then deleting small open neighbourhoods of its singular vertices yields a compact $n$-manifold $\check N$ with no spherical boundary components.
Conversely, from any compact $n$-manifold $M$, a singular $n$-manifold $\widehat M$ can be constructed by capping off each component of $\partial M$ by a cone over it.
\\
Moreover, within the class of compact $n$-manifolds with no spherical boundary components, that we will consider throughout the present work,
the above correspondence is well-defined, therefore singular $n$-manifolds and compact $n$-manifolds of this class can be bijectively associated.}
\end{rem}

\begin{defn}\label{colored triangulation}
{\em Given a singular $n$-manifold $N$, a {\em colored triangulation} of $N$ is a pseudocomplex\footnote{The notion of pseudocomplex is a generalization of that of simplicial complex.
Roughly speaking, in a pseudocomplex two $n$-simplexes may have more than one common $(n-1)$-face; however identifications of faces of the same simplex are not allowed.} 
$K$ triangulating $N$ and carrying a vertex-labelling by means of the elements of the  set $\Delta_n=\{0,\ldots,n\},$ such that vertices belonging to the same simplex have different labels.
\\ \noindent Given an $h$-simplex ($0\leq h\leq n$) $\sigma$ of $K$, the {\em disjoint star} of $\sigma$ in $K$ is the pseudocomplex obtained by taking the $n$-simplexes of $K$ 
which contain $\sigma$ and identifying, according to the gluings in $K$, only their $(n-1)$-faces containing $\sigma.$ The {\em disjoint link} of $\sigma$ in $K$ is the subcomplex of its disjoint star
formed by the simplexes that do not contain $\sigma.$}
\end{defn}

\begin{defn} \label{$n+1$-colored graph}
{\em An {\em $(n+1)$-colored graph}  ($n \ge 2$) is a pair $(\Gamma,\gamma)$, where $\Gamma=(V(\Gamma), E(\Gamma))$ is an $(n+1)$-regular multigraph (i.e. multiple edges are allowed, but loops are forbidden) 
and $\gamma$ is an {\em edge-coloration}, that is a map  $\gamma: E(\Gamma) \rightarrow \Delta_n=\{0,\ldots,n\}$ assigning different colors to incident edges. %PC anche to edges sharing an endpoint
}\end{defn}

Given $\{c_1,\dots,c_h\}\subseteq\Delta_n$, let us set $\hat c_1\cdots\hat c_h = \Delta_n\setminus\{c_1,\dots,c_h\}$; by $\Gamma_{c_1,\dots,c_h}$ 
(resp. $\Gamma_{\hat c_1\cdots\hat c_h}$) we will denote the spanning subgraph of $(\Gamma, \gamma)$ whose edge set contains only the edges of $\G$ that are colored  
(resp. that are not colored) by $\{c_1,\dots,c_h\}$. 
The connected components of $\Gamma_{c_1,\dots,c_h}$ will be called {\it $\{c_1,\dots,c_h\}$-residues} %or {\it $h$-residues} 
of $\G$ and their number will be denoted by 
$g_{c_1,\dots,c_h}$. Furthermore, we will use the notation $g_{\hat c_1\cdots\hat c_h}$ for the number of connected components of $\Gamma_{\hat c_1\cdots\hat c_h}.$
\medskip

Given a colored triangulation $K$ of a singular $n$-manifold $N$, the dual 1-skeleton of $K$ is an $(n+1)$-colored graph $(\G(K),\g(K))$, whose edges inherit the coloration from 
the vertex-labelling of $K;$ more precisely, for each $c\in\Delta_n,$ an edge $e\in E(\G(K))$ is $c$-colored iff %the vertices of 
the $(n-1)$-simplex of $K$ dual to $e$ misses color $c.$

We will say that $(\G(K),\g(K))$ {\it represents} the singular manifold $N$ and, in virtue of Remark \ref{correspondence-sing-boundary}, also the compact $n$-manifold $\check N.$

The duality between $K$ and $(\G(K),\g(K))$ gives rise to a bijective correspondence between $(n-h)$-simplexes of $K$ and $h$-residues of $\G(K);$
more precisely, to each $(n-h)$-simplex $\sigma$ of $K$ whose vertices are labelled by $\hat c_1\cdots\hat c_{h}$ corresponds the $\{c_1,\dots,c_{h}\}$-residue of $\G(K)$ 
representing the disjoint link of $\sigma$ in $K.$
Conversely, any $(n+1)$-colored graph $(\Gamma,\gamma)$ can be thought as the dual 1-skeleton of a pseudocomplex $K(\G)$, whose $n$-simplexes are in bijection with the
vertices of $\G$ and the gluings of $(n-1)$-faces are determined by the colored adjacencies of $\G$ (see \cite{Ferri-Gagliardi-Grasselli} for a more detailed description of this construction).

By the duality between $(\Gamma,\gamma)$ and $K(\G)$ it is not difficult to see that $|K(\G)|$ is a singular $n$-manifold iff for any $c\in\Delta_n$, each $\hat c$-residue
of $\G$ represents a closed connected $(n-1)$-manifold\footnote{Note that, in this case, each $\{c_1,\dots,c_h\}$-residue ($1\leq h < n$) of $\G$ represents an $(n-h-1)$-sphere.};
residues not representing $(n-1)$-spheres are called {\it singular}.
In particular, $|K(\G)|$ is a closed $n$-manifold iff for any $c\in\Delta_n$, no $\hat c$-residue of $\G$ is singular.
\bigskip

Throughout the paper we will be interested in colored graphs representing singular $n$-manifolds with at most one singular vertex, whose associated compact 
manifolds have therefore empty or connected boundary. In particular, we will consider $(n+1)$-colored graphs whose dual colored triangulations have the minimum possible number
of vertices.

\begin{defn} {\em An $(n+1)$-colored graph $(\G,\g)$ representing a compact $n$-manifold $M$ with empty or connected boundary, is called a {\em crystallization} of $M$ if
$g_{\hat c}=1,\ \ \forall c\in\Delta_n.$ 
\\ \noindent Its dual colored triangulation $K(\G)$, having exactly $n+1$ vertices, is said to be {\em contracted}.} 
\end{defn}

\begin{rem}\label{rem-crystallization} {\em The existence of combinatorial moves ({\em dipole moves}) on edge-colored graphs inducing PL homeomorphisms between the represented compact manifolds,
allows to prove that {\em any compact $n$-manifold with empty or connected boundary admits a crystallization} (\cite{Casali-Cristofori-Gagliardi_Rendiconti_TS_2020}).} 
\end{rem}

The possibility of defining, for compact manifolds, the PL invariant regular genus relies on the existence of a particular type of embedding of colored graphs into surfaces. 

\begin{defn} {\em A cellular embedding of an $(n+1)$-colored graph $(\Gamma,\gamma)$ into a closed surface is called {\it regular} if there exists a cyclic permutation $\varepsilon$ 
of $\Delta_n$ such that the regions of the embedding are bounded by the images of the  cycles of $\G$ whose edges are colored by two colors that are consecutive in $\varepsilon.$
}\end{defn}
 
An important result by Gagliardi (\cite{Gagliardi 1981}) guarantees that given a bipartite (resp. non-bipartite) $(n+1)$-colored graph $(\Gamma,\gamma)$, for each cyclic 
permutation $\varepsilon = (\varepsilon_0,\ldots,\varepsilon_n)$ of $\Delta_n$ there exists a regular embedding of $(\Gamma,\gamma)$ such that the embedding surface
$F_{\varepsilon}(\Gamma)$ is orientable (resp. non-orientable) and has genus (resp. half genus) $\rho_{\varepsilon} (\Gamma)$ satisfying 

\begin{equation*}
2 - 2\rho_\varepsilon(\Gamma)= \sum_{j\in \mathbb{Z}_{n+1}} g_{\varepsilon_j,\varepsilon_{j+1}} + (1-n)p,
\end{equation*}

\noindent where $2p$ is the order of $\G.$  

\noindent $\rho_{\varepsilon} (\Gamma)$ is called the {\it regular genus of $\Gamma$ with respect to the permutation $\varepsilon$}.

Moreover, $\varepsilon$ and $\varepsilon^{-1}$ induce the same embedding and $(\Gamma,\gamma)$ admits no regular embeddings into non-orientable (resp. orientable) surfaces.

\bigskip

As a consequence, the following definitions can be given:

\begin{defn} \label{regular-genus}
{\em The {\em regular genus} $\rho(\Gamma)$ of an $(n+1)$-colored graph $(\Gamma,\gamma)$ is defined as}
\begin{equation*}
 \rho(\Gamma) \ = \ \min\, \{\rho_\varepsilon(\Gamma)\ /\ \varepsilon \mbox{\ is a cyclic permutation of\ } \Delta_n\}.
\end{equation*}
{\em The {\em (generalized) regular genus} of a compact $n$-manifold $M$ ($n\geq 2$) is defined as}   
\begin{equation*}
\mathcal G(M)=\min \{\rho(\Gamma)\ | \ (\Gamma,\gamma)\mbox{ represents} \ M\}.
\end{equation*}
\end{defn}

\bigskip

%\begin{definition}\label{Gdegree} 
%Let $(\Gamma,\gamma)$ be an $(n+1)$-colored graph and let $\varepsilon^{(1)}, \varepsilon^{(2)}, \dots , \varepsilon^{(\frac {n!} 2)}$ be the cyclic permutations of $\Delta_n$ 
%up to inverse. 
%The \emph{Gurau degree} (or \emph{G-degree} for short) $\omega_{G}(\Gamma)$ of $\Gamma$ is defined as
%\begin{equation*}
% \omega_{G}(\Gamma) \ = \ \sum_{i=1}^{\frac {n!} 2} \rho_{\varepsilon^{(i)}}(\Gamma)
%\end{equation*}
%Let $M$ be a compact (PL) $n$-manifold ($n\geq 2$). The {\it Gurau degree} (or {\it G-degree}) of $M$ is defined as
%\begin{equation*}
%\mathcal D_G(M)=\min \{\omega_G(\Gamma)\ | \ (\Gamma,\gamma)\mbox{ represents} \ M\}.
%\end{equation*}
%\end{definition}

\begin{rem}{\em Any bipartite (resp. non-bipartite) $3$-colored graph $(\Gamma,\gamma)$ represents an orientable (resp. non-orientable) surface $|K(\Gamma)|$ and 
$\rho(\Gamma)$ coincides with the genus (resp. half the genus) of $|K(\Gamma)|;$ 
moreover, the regular genus of a closed 3-manifold coincides with its Heegaard genus and the coincidence still holds in the case of connected boundary (\cite{Gagliardi 1981}, 
\cite{Cristofori-Gagliardi-Grasselli}) making the regular genus an %interesting/valid/convincing/
effective extension to dimension $n>3$ of these classical notions of genus. 
\\
Actually many classification results according to this invariant already exist in dimension $4$ and $5$ (\cite{Casali_Forum2003}, \cite{Casali-Gagliardi ProcAMS}, 
\cite{Casali-Cristofori-Gagliardi Complutense 2015}, \cite{Casali-Cristofori_generalized-genus} and references therein).
}
\end{rem}

\medskip

Presentations of the fundamental groups of a compact $n$-manifold and its associated singular manifold can be obtained from any
$(n+1)$-colored graph representing them in the following way (\cite{Casali-Cristofori_generalized-genus}): 

\begin{prop}  \label{fundamental_group} 
Let $(\Gamma,\gamma)$ be an $(n+1)$-colored graph representing a compact $n$-manifold $M$. 

\noindent For each $i,j \in \Delta_n$,  let $X_{ij}$  (resp. $R_{ij}$) be a set in bijection with the connected components of $\Gamma_{\hat i \hat j}$ (resp. with the $\{i,j\}$-colored 
cycles of $\Gamma$), and let $\bar R_{ij}$ be a subset  of $X_{ij}$ corresponding to the a maximal tree of the subcomplex $K_{ij}$ of $K(\Gamma)$ (consisting only of vertices 
labelled $i$ and $j$, and edges connecting them). 
Then: 
\begin{itemize}
\item[(a)]
if $i,j\in \Delta_n$ are not singular in $\Gamma$, \ \ \ $$  \pi_1(M) \ = \ < X_{ij} \ / \ R_{ij} \cup \bar R_{ij}>;$$ 
\item[(a')]
if no color in $\Delta_n - \{i,j\}$ is singular in $\Gamma$,  \ \ \ $$  \pi_1(\widehat M) \ = \ < X_{ij} \ / \ R_{ij} \cup \bar R_{ij} >.$$
\end{itemize}
\end{prop}

\bigskip
\bigskip

\section{Weak simple crystallizations of compact $4$-manifolds}\label{sec-weaksimple}

In \cite{Basak-Casali} and \cite{Basak} {\it semi-simple} and {\it weak semi-simple} crystallizations of closed $4$-manifolds are introduced and studied, 
by generalizing the former notion of {\it simple crystallizations} of closed simply-connected $4$-manifolds (see \cite{Basak-Spreer} 
and \cite{Casali-Cristofori-GagliardiJKTR2015}): they are proved to be ``minimal" with respect to regular genus, among all graphs representing the same closed $4$-manifold.
In  \cite{Casali-Cristofori_trisections} these definitions are extended to compact $4$-manifolds with empty or connected boundary and their minimizing property with 
respect to regular genus is proved to hold also in this wider context (\cite{Casali-Cristofori-Gagliardi_Rendiconti_TS_2020}).\footnote{Actually, weak semi-simple crystallizations  turn out to admit also interesting properties with respect to the existence of {\it trisections} of the represented compact 4-manifolds, and to the computation of their {\it trisection genus}: see \cite{Casali-Cristofori_trisections} and references therein.}    %MR_26agosto (anche se mi sono accorta che l'articolo sulle trisezioni era comunque già citato!) 

Since the present paper focuses on simply-connected $4$-manifolds (with empty or connected boundary), we will restrict the definitions to this particular case, where the above
types of crystallizations %involved objects 
are called {\it simple} and {\it weak simple} %crystallizations 
respectively.

%In the particular case of simply-connected $4$-manifolds (with empty or connected boundary), they are called {\it simple} and {\it weak simple} crystallizations respectively. %, and their definitions are as follows: 

\bigskip
\bigskip

%From now on and till the end of this section we will restrict our attention to compact simply-connected $4$-manifolds with empty or connected  boundary.  
In the following, given a $5$-colored graph $(\Gamma,\gamma)$ representing a compact $4$-manifold with empty or connected  boundary, we %may 
will always assume, without loss of generality, %(see Remark \ref{rem-crystallization}),  we may assume $\Gamma$  to be a crystallization %(i.e. $\Gamma_{\hat c}$ is connected for any $c \in \Delta_4$) and
that color $4$ is its (unique) possible singular color. %(i.e. $\Gamma_{\hat c}$ represents $\mathbb S^3$, for any $c \ne 4$).  

Furthermore, $\mathcal P_4$ will denote the set of all cyclic permutations of $\Delta_4$, and, given $\varepsilon=(\varepsilon_0,\varepsilon_1,\varepsilon_2,\varepsilon_3,\varepsilon_4)\in\mathcal P_4$, 
we will always suppose, again without loss of generality, $\varepsilon_4=4.$ 

\medskip

We can now introduce the following definitions:

\begin{defn} \label{def_weak simple} {\rm Let $M^4$ be a compact simply-connected $4$-manifold, with empty or connected boundary.  A $5$-colored graph $(\Gamma,\gamma)$ representing $M^4$   
 is called {\em weak simple} with respect to $\varepsilon \in \mathcal P_4$  %a permutation
if \, $g_{\varepsilon_i,\varepsilon_{i+2},\varepsilon_{i+4}} = 1 \ \ \forall \ i \in \Delta_4$ (where the additions in subscripts are intended in $\mathbb{Z}_{5}$).
$(\Gamma,\gamma)$ is called {\em simple} if \, $g_{j,k,l} = 1\ \ \forall \ j,k,l \in \Delta_4.$
}
\end{defn}
\smallskip

Note that, as a consequence of the above definition, if $(\Gamma,\gamma)$ is weak simple, then \ $g_{\hat j}=1,\ \forall \ j\in\Delta_4$, i.e. $(\Gamma,\gamma)$ is a crystallization of $M^4.$

\medskip

Let $(\Gamma,\gamma)$ be a crystallization of a compact simply-connected $4$-manifold $M^4$ with empty or connected boundary (i.e. $\Gamma_{\hat c}$ represents $\mathbb S^3$, for any $c \ne 4$),
and let us denote, for shortness, by $\rho_\e$  (resp. by $\rho_{\e_{\hat{i}}}$) the regular genus of $\G$ (resp. of $\ \G_{\widehat{\e_i}}$)
with respect to the permutation $\varepsilon$  (resp. $\ \e_{\hat i} = (\e_0,\ldots,\e_{i-1},\e_{i+1},\ldots,\e_4=4)$).

By setting $g_{j,k,l} = 1 + t_{j,k,l}\quad j,k,l\in\Delta_4$ (see \cite{Basak} and \cite{Casali-Cristofori_trisections}), the following equality can be proved (see \cite{Casali-Cristofori_trisections} and 
\cite{Casali-Cristofori-Gagliardi_Rendiconti_TS_2020}):

\begin{equation} \label{genus-subgenera(t)}
\rho_{\varepsilon}  - \rho_{\varepsilon_{\hat i}}  -  \rho_{\varepsilon_{\widehat {i+2}}} =  \ t_{\varepsilon_{i-1},\varepsilon_{i+1},\varepsilon_{{i+3}}}\geq 0  \ \  \ \   \forall i \in \Delta_4
\end{equation}

\noindent where all subscripts are taken in $\mathbb Z_5.$  

\bigskip

\bigskip

By making use of relation  \eqref{genus-subgenera(t)}, it is not difficult to prove the following characterization of weak simple crystallizations: 

\begin{prop} {\rm (\cite[Corollary 8]{Casali-Cristofori_trisections})}\label{weaksimple}
Let $(\Gamma,\gamma)$ be a crystallization of a compact simply-connected $4$-manifold $M^4$ with empty or connected boundary. 
Then $(\Gamma,\gamma)$ is weak simple with respect to %a cyclic permutation 
$\varepsilon \in \mathcal P_4$ if and only if for each  $i\in  \Delta_4,\ \rho_{\varepsilon_{\hat i}} = \frac 1 2 \rho_{\varepsilon}.$
\end{prop}

\bigskip
\bigskip
%The Main Theorem of \cite{Casali-Cristofori-Gagliardi_Rendiconti_TS_2020} explains the minimizing properties of semi-simple crystallizations (resp. of weak semi-simple crystallizations) with respect to the PL invariants regular genus, G-degree and gem-complexity (resp. with respect to the PL invariants regular genus).
The Main Theorem of \cite{Casali-Cristofori-Gagliardi_Rendiconti_TS_2020} explains the minimizing property of weak semi-simple crystallizations with respect to regular genus.
In the simply-connected case, the result becomes:

\begin{thm} {\rm \cite[Main Theorem]{Casali-Cristofori-Gagliardi_Rendiconti_TS_2020}} \label{upper_bound_total}
Let $M^4$ be a compact simply-connected $4$-manifold  with empty or connected boundary, %and let $\widehat M^4$ be its associated singular manifold. 
then the regular genus $ \mathcal G (M^4)$ of $M^4$ satisfies     
$$ \mathcal G (M^4) \ge 2 \chi(\widehat M^4) -4.$$
Moreover, equality holds if and only if $M^4$ admits a weak simple crystallization. 
\end{thm}

\bigskip

\begin{rem}{\em We point out that simple, or semi-simple in the not simply-connected case, crystallizations realize also the PL invariants {\em gem-complexity} and {\em G-degree} (\cite{Casali-Cristofori-Gagliardi_Rendiconti_TS_2020}).
Gem-complexity of a compact PL $4$-manifold $M^4$ is related to the minimum number of vertices in a $5$-colored graph representing $M^4$, while the G-degree is a PL invariant arising
within the theory of {\em Colored Tensor Models}, that was first developed, in theoretical physics, as a random geometry approach to quantum gravity and has subsequently gained interest
in the study of quantum mechanical models, too (see \cite{Gurau-book}, \cite{Casali-Cristofori-Dartois-Grasselli}, \cite{Gurau2017}, \cite{Witten}).} 
\end{rem}

\bigskip
\medskip

As a consequence of Theorem \ref{upper_bound_total}, together with Proposition \ref{weaksimple}, we have: 

\begin{cor}\label{weaksimple minimal genus}
Let $(\Gamma,\gamma)$ be a crystallization of % $5$-colored graph representing 
a compact simply-connected $4$-manifold $M^4$ with empty or connected boundary, then, $(\Gamma,\gamma)$ is weak simple with respect to %the cyclic permutation 
$\varepsilon \in \mathcal P_4$ if and only if             
$$ \mathcal G (M^4) = \rho (\Gamma) = \rho_{\varepsilon}(\Gamma) = 2 \chi(\widehat M^4) -4 $$
or, equivalently, if and only if 
$\ \rho_{\varepsilon_{\hat i}} = \chi(\widehat M^4) -2 \ \ \forall i \in  \Delta_4.$
\end{cor}
\ \qed

\bigskip

\begin{prop}\label{reg-gen-Betti-number} If $M^4$ is a simply-connected $4$-manifold with empty or connected boundary, then 
$$\mathcal G (M^4) \ge 2\beta_2(M^4),$$ where $\beta_2(M^4)$ is the second Betti number of $M^4.$ 

% Moreover, if $\G$ is a crystallization of $M^4$, then $\rho(\Gamma) = \rho_\varepsilon(\Gamma) = 2\beta_2(M^4)$ iff $\Gamma$ is weak simple with respect to a cyclic permutation $\varepsilon$ of $\mathcal P_4$. %of $M^4$ if and only if equality holds.  In this case all $\hat\varepsilon_i$-residues ($i\in\mathbb Z_5$) of $\G$ have regular genus $\beta_2(M^4)$ with respect to the permutations induced by  $\varepsilon.$
Moreover, $\mathcal G (M^4) = 2\beta_2(M^4)$ holds if and only if $M^4$ admits a  weak simple crystallization $(\Gamma,\gamma)$ with respect to %a cyclic permutation 
$\varepsilon$ of $\mathcal P_4$ or, equivalently, 
if and only if all $\hat\varepsilon_i$-residues ($i\in\mathbb Z_5$) of $\G$ have regular genus $\beta_2(M^4)$ with respect to the permutations induced by  $\varepsilon.$
\end{prop}

\dimo
First, let us observe that, for any compact $4$-manifold $M^4$, if $C$ is the cone over $\partial M^4$, then the inclusion $(M,\partial M)\hookrightarrow (\widehat M^4,C)$ induces an isomorphism in homology;
therefore the exact sequence of relative homology of the pair $(\widehat M^4,C)$, yields 
$$H_k(\widehat M^4)\cong H_k(\widehat M^4,C)\cong H_k(M^4,\partial M^4)\cong H^{4-k}(M^4)\qquad \forall k\in\{2,3\}$$
\noindent where the last isomorphism comes from Poincar\'e-Lefschetz duality.

Hence $\beta_2(\widehat M^4)=\beta_2(M^4),\ \beta_3(\widehat M^4)=\beta_1(M^4)$ and, therefore 
\begin{equation} \label{Euler characteristic}
\chi(\widehat M^4) = 2  - \beta_1(\widehat M^4) + \beta_2(M^4) - \beta_1(M^4).
\end{equation}  

By further assuming $M^4$ to be simply-connected with empty or connected boundary, statement (a) of the Main Theorem implies 
$$\mathcal G (M^4) \ge 2\chi(\widehat M^4)- 4 = 2\beta_2(M^4).$$

Now, %if $\Gamma$ is a crystallization of $M^4$, then by Corollary \ref{weaksimple minimal genus}, %Proposition \ref{invariants computation}(a), 
by Corollary \ref{weaksimple minimal genus},  a crystallization $\Gamma$ of $M^4$ realizes the lower bound for regular genus, i.e. 
$\rho(\Gamma) = 2\beta_2(M^4)$, 
if and only if there exists a cyclic permutation $\varepsilon$ of $\mathcal P_4$ such that $\Gamma$ is weak simple with respect to $\varepsilon$ 
(in this case $\rho(\Gamma) = \rho_\varepsilon(\Gamma)$).

%With this assumption, 
Finally, Corollary \ref{weaksimple minimal genus} ensures the equivalence with the condition
$$\rho_{\varepsilon_{\hat i}} = \chi(\widehat M^4) - 2 = \beta_2(M^4)\ \ \forall i \in  \Delta_4.$$
\  \qed 

\bigskip

\section{Proof of the main results}

In order to prove the main results of the paper, we need two lemmas.

\begin{lemma} \label{lemma-beta_2}
Let $M^4$ be a compact simply-connected 4-manifold. Then, for each 5-colored graph $(\Gamma,\gamma)$ representing $M^4$ and for each cyclic permutation $\e$ of $\Delta_4$,
$$ \sum_{i \in \Delta_4} \rho_{\varepsilon_{\hat i}}(\Gamma_{\widehat{ \varepsilon_i}}) - 2 \rho_\e(\Gamma) =  \beta_2(M^4).$$  
As a consequence, 
$$ \beta_2(M^4) \le min  \{ \rho_{\varepsilon_{\hat i}}(\Gamma_{\widehat{ \varepsilon_i}}) \ / \ (\Gamma,\gamma) \ \text{representing} \ M^4, \ \e \ \text{cyclic permutation of}  \ \Delta_4, \ i \in \Delta_4 \}.$$
\end{lemma}

\dimo
According to formula (20) of Proposition 13 in \cite{Casali-Cristofori_generalized-genus}, for each 5-colored graph representing a compact 4-manifold $M^4$ and for each cyclic permutation $\e$ of $\Delta_4$, the following relation holds: 
\begin{equation} \label{Euler-caracteristic n=4}
 %\chi (\widehat M^4)  = 2 - 2 \rho_{\varepsilon}(\Gamma) + \sum_{i \in \Delta_4} \rho_{\varepsilon}(\Gamma_{\widehat{ \varepsilon_i}}).
 \chi (\widehat M^4)  = 2 - 2 \rho_{\varepsilon} + \sum_{i \in \Delta_4} \rho_{\varepsilon_{\hat i}}.
\end{equation}

By comparing with formula \eqref{Euler characteristic}, we obtain: 
\begin{equation} \label{genus & sum_subgenus vs Betti_numbers}
%\sum_{i \in \Delta_4} \rho_{\varepsilon}(\Gamma_{\widehat{ \varepsilon_i}}) - 2 \rho_\e(\Gamma) =  \beta_2(M^4) - \beta_1(M^4) - \beta_1(\widehat M^4).
\sum_{i \in \Delta_4} \rho_{\varepsilon_{\hat i}} - 2 \rho_\e =  \beta_2(M^4) - \beta_1(M^4) - \beta_1(\widehat M^4).
\end{equation}
If $M^4$ is assumed to be simply-connected, the general relation of the statement directly follows.
% $$ \sum_{i \in \Delta_4} \rho_{\varepsilon}(\Gamma_{\widehat{ \varepsilon_i}}) - 2 \rho_\e(\Gamma) =  \beta_2(M^4).$$   

On the other hand, it is well known that, for each 5-colored graph representing a compact 4-manifold $M^4$ and for each $\e=(\e_0,\e_1,\e_2,\e_3,\e_4)\in\mathcal P_4,$
the inequality 
$$  \rho_{\e_{\widehat {j-1}}}  +  \rho_{\e_{\widehat {j+1}}} \le \rho_{\e} $$  
holds $\forall j \in \Delta_4$ %, where $\rho_{\e_{\hat i}}$ denotes $\rho_{\varepsilon}(\Gamma_{\widehat{ \varepsilon_i}})$ ($\forall i\in \Delta_4$)
: see  formula (24) of Proposition 7.1 in  \cite{Casali-Cristofori_generalized-genus}.

Hence, the hypothesis of simply-connectedness implies, $\forall i \in \Delta_4$: 
%$$ \beta_2 (M^4)=  \sum_{i \in \Delta_4} \rho_{\varepsilon}(\Gamma_{\widehat{ \varepsilon_i}}) - 2 \rho_\e(\Gamma)   =   \rho_{\e_{\hat i}}  + ( \rho_{\e_{\widehat i-1}}  +  \rho_{\e_{\widehat {i+2}}}) + ( \rho_{\e_{\widehat i-2}}  +  \rho_{\e_{\widehat {i+1}}}) - 2 \rho_\e(\Gamma) \le \rho_{\e_{\hat i}}.$$
$$ \beta_2 (M^4)=  \sum_{i \in \Delta_4} \rho_{\varepsilon_{\hat i}} - 2 \rho_\e   =   \rho_{\e_{\hat i}}  + ( \rho_{\e_{\widehat {i-1}}}  +  \rho_{\e_{\widehat {i+2}}}) + ( \rho_{\e_{\widehat {i-2}}}  +  \rho_{\e_{\widehat {i+1}}}) - 2 \rho_\e \le \rho_{\e_{\hat i}}.$$
\ \ \qed

% \begin{lemma} \label{lemma_subgenus} 
% Let $(\Gamma, \gamma)$ be a crystallization of a (simply-connected) compact 4-manifold $M^4$ (with empty or  connected boundary), such that color $4$ is the only possible singular color. 
% If there exist  two different colors $j,k\in \Delta_3$ such that $ g_{\hat j \hat k}= 1$, then  for each $s \in \Delta_3-\{j,k\},$ a cyclic permutation $\e$ of $\Delta_4$ exists, with $\e_0 = s$, so that: 
% $$ \rho_{\e_{\hat 0}} = \beta_2(M^4) + t_{s,j,k}.$$

% In particular, if there exist  three different colors $r,j,k\in \Delta_3$ such that $g_{\hat r \hat 4}=  g_{\hat j \hat k}= 1$, then a cyclic permutation $\e$ of $\Delta_4$ exists, so that $ \rho_{\e_{\hat 0}} = \beta_2(M^4)$, where $\{\e_0\}= \Delta_4 - \{r,j,k,4\}.$ 
%\end{lemma}

\begin{lemma} \label{lemma_subgenus}
Let $(\Gamma, \gamma)$ be a crystallization of a compact 4-manifold $M^4$ (with empty or  connected boundary), such that color $4$ is the only possible singular color. 
If there exist  two different colors $j,k\in \Delta_3$ such that $ g_{\hat j \hat k}= 1$, then $M^4$ is simply-connected and for each $s \in \Delta_3-\{j,k\},$ a cyclic permutation $\e$ of $\Delta_4$ exists, with $\e_0 = s$, so that: 
$$ \rho_{\e_{\hat 0}} = \beta_2(M^4) + t_{s,j,k}.$$

In particular, if there exist  three different colors $r,j,k\in \Delta_3$ such that $g_{\hat r \hat 4}=  g_{\hat j \hat k}= 1$, then a cyclic permutation $\e$ of $\Delta_4$ exists, so that $ \rho_{\e_{\hat 0}} = \beta_2(M^4)$, where $\{\e_0\}= \Delta_4 - \{r,j,k,4\}.$ 
\end{lemma}

\dimo
According to Section 4 of  \cite{Casali-Cristofori_trisections}, for each crystallization of a compact 4-manifold with either empty or connected boundary, setting $m = rk(\pi_1(M^4))$ and $m ^\prime = rk(\pi_1(\widehat{M^4}))$, we have: 
$$ g_{\e_{i-1},\e_{i+1},\e_{i+3}} = 1 + \rho_\e  - \rho_{\e_{\hat i}}   -  \rho_{\e_{\widehat{i+2}}}  =  1+m^{\prime} + t_{\e_{i-1},\e_{i+1},\e_{i+3}}  \ \  \ \   \forall i \in \{2,4\} \ \ \text{and} $$
$$ g_{\e_{i-1},\e_{i+1},\e_{i+3}} = 1 + \rho_\e  - \rho_{\e_{\hat i}}   -  \rho_{\e_{\widehat{i+2}}}  =  1+m + t_{\e_{i-1},\e_{i+1},\e_{i+3}} \ \  \ \   \forall i \in \{0,1,3\},$$
which trivially imply 
% \begin{equation} %\label{genus-subgenera}
$$\begin{aligned} \rho_{\e}  - \rho_{\e_{\hat i}}  -  \rho_{\e_{\widehat {i+2}}} -  m^{\prime} = & \ t_{\e_{i-1},\e_{i+1},\e_{{i+3}}}  \ \  \ \   \forall i \in \{2,4\}  \ \ \text{and} \\
 \rho_{\e}  - \rho_{\e_{\hat i}}  -  \rho_{\e_{\widehat {i+2}}} -  m =  & \ t_{\e_{i-1},\e_{i+1},\e_{{i+3}}}   \ \  \ \   \forall i \in \{0,1,3\}. \end{aligned}$$ % \end{equation} 
  
%Hence, if $\e$ denotes a cyclic permutation of $\Delta_4$ such that both $r$ and $4$ and  $j,k$ are not consecutive, $r$ being the only element of $\Delta_2-\{s,j,k\}$ (for example: 
%$\e=(\e_0=s,\e_1=j,\e_2=r,\e_3=k,\e_4=4)$), 
% then $g_{\hat j \hat k}= 1$ implies $m=0$, $t_{s,r,4} =0$ and $\rho_{\hat j} + \rho_{\hat k} =\rho_{\e}.$
%On the other hand, $m=0$ implies $m^\prime=0$, and hence $\rho_{\hat r} + \rho_{\hat 4} =\rho_{\e}- t_{s,j,k}  $
%As a consequence, $M^4$ is simply-connected and Lemma \ref{lemma-beta_2}  easily yields \ $  \rho_{\hat s}  - t_{s,j,k} = \beta_2(M^4).$  

Let $\e$ be a cyclic permutation of $\Delta_4$ with $\e_0 = s$ and such that both $r,4$ and  
$j,k$ are not consecutive, $r$ being the only element of $\Delta_2-\{s,j,k\}$; for example, let us suppose : $\e=(\e_0=s,\e_1=j,\e_2=r,\e_3=k,\e_4=4)$), 
then $g_{\hat j \hat k}= 1$ implies $m=0$ (see Proposition \ref{fundamental_group}(a)), $t_{s,r,4} =0$ and $\rho_{\e_{\hat 1}} + \rho_{\e_{\hat 3}} =\rho_{\e}.$   
On the other hand, $m=0$ implies $m^\prime=0$, and hence $\rho_{\e_{\hat 2}} + \rho_{\e_{\hat 4}} =\rho_{\e}- t_{s,j,k}  $    
As a consequence, $M^4$ is simply-connected and Lemma \ref{lemma-beta_2}  easily yields \ $  \rho_{\e_{\hat 0}}  - t_{s,j,k} = \beta_2(M^4).$ 

The second part of the statement directly follows from the first one, by recalling that   $g_{\hat r \hat 4}=1 $ means %\  \Leftrightarrow \ 
$t_{s,j,k}=0$. \ \ \qed 

\bigskip

Let us recall now that that every closed (resp. compact) 4-manifold $M^4$ admits a handle decomposition
$$ M^4= H^{(0)} \cup (H^{(1)}_1 \cup \dots \cup H^{(1)}_{r_1}) \cup (H^{(2)}_1 \cup \dots \cup H^{(2)}_{r_2}) \cup (H^{(3)}_1 \cup \dots \cup H^{(3)}_{r_3}) \cup  H^{(4)}$$
(resp. 

$$ M^4= H^{(0)} \cup (H^{(1)}_1 \cup \dots \cup H^{(1)}_{r_1}) \cup (H^{(2)}_1 \cup \dots \cup H^{(2)}_{r_2}) \cup (H^{(3)}_1 \cup \dots \cup H^{(3)}_{r_3})) $$

where $H^{(0)} =\mathbb D^4$ and  each $p$-handle $H^{(p)}_i= \mathbb D^p \times \mathbb D^{4-p}$  ($1 \le p \le 4, \ 1 \le i \le r_p$) is endowed with an embedding (called {\it attaching map}) $f_i^{(p)}: \partial \mathbb D^p \times \mathbb D^{4-p} \to  \partial (H^{(0)}  \cup \dots (H^{(p-1)}_1 \cup \dots \cup H^{(p-1)}_{r_{p-1}}))$; 
moreover, in the closed case, it is well-known that $3$- and $4$-handles are attached in a unique way to the union of the $h$-handles, with $0 \leq h \leq 2$.

A standard argument of crystallization theory allows to state that, for any crystallization $(\Gamma,\gamma)$ of a closed 4-manifold $M^4$ and for any partition $\{\{i,j,k\},\{r,s\}\}$ of $\Delta_4$,  $M^4$ admits a decomposition of type $ M^4 \, = \, N(i,j,k) \cup_{\phi} N(r,s),$  where $N(i,j,k)$ 
(resp. $N(r,s)$) denotes a regular neighbourhood of the subcomplex $K(i,j,k)$ (resp. $K(r,s)$) of $K(\Gamma)$ generated by the vertices labelled $\{i,j,k\}$ (resp. $\{r,s\}$) 
and $\phi$ is a boundary identification (see \cite{Ferri-Gagliardi Proc AMS 1982} and \cite{Gagliardi 1981}).
The same decomposition applies to the singular manifold $\widehat M^4,$ in case $(\Gamma,\gamma)$ is a crystallization of a compact 4-manifold. 

Any such decomposition turns out to induce a handle decomposition of the closed (resp. compact, with non-empty boundary) $4$-manifold $M^4$, where $N(i,j,k)$ constitutes the 
union of the $h$-handles, with $0 \leq h \leq 2$, while $N(r,s)$ is the union of $3$- and $4$-handles (resp. is the union of $3$-handles together with the cones over the
components of $\partial M^4$).
% In the particular case of a simple crystallization, it is easy to prove that a so-called {\it special handlebody decomposition} of $M$ is actually obtained, i.e. a handle decomposition lacking in 1-handles and 3-handles (see \cite[Section~3.3]{[M]}):
\medskip

The following proposition establishes combinatorial conditions on $(\Gamma,\gamma)$ which ensure the existence of particular types of handle decompositions of the represented manifold (lacking in 1-handles).

\begin{prop} \label{handle-decomposition_general}
Let $(\Gamma, \gamma)$ be a crystallization of a (simply-connected) compact 4-manifold $M^4$ (with empty or  connected boundary), such that color $4$ is the only possible singular color. 
If there exist three different colors $i,j,k\in \Delta_3$ such that $g_{\hat i \hat k} =g_{\hat j \hat k}= 1$, then  $M^4$ admits a handle decomposition consisting of one 0-handle, $\beta_2(M^4)+t_{i,j,k}$ 2-handles, $t_{i,j,k}$ 3-handles and - in the closed case - one 4-handle.

In particular, if there exists a partition $\{\{4,r\}, \{i,j,k\}\}$ of $\Delta_4$ such that $g_{\hat r \hat 4}=  g_{\hat i \hat k} =g_{\hat j \hat k}= 1$, then the induced handle decomposition of $M^4$ consists of one 0-handle, $\beta_2(M^4)$ 2-handles and - in the closed case - one 4-handle.
\end{prop}

\dimo
Let us denote by $\varepsilon$ the following cyclic permutation of $\Delta_4$:   $\e=(\e_0=i,\e_1=j,\e_2=r,\e_3=k,\e_4=4),$ where $r$ is the only element of $\Delta_3-\{i,j,k\}$.  
Then, let us take into account % the partition $\{\{4,r\}, \{i,j,k\}\}$ of $\Delta_4$ coincides with 
the partition $\{\{\e_0,\e_1,\e_3\},\{\e_2,\e_4\}\}$ of $\Delta_4$ (i.e. $\{\{4,r\}, \{i,j,k\}\}$).  
Now, in order to analyze the associated handle decomposition, first of all note that $K(\e_2,\e_4)=K(r,4)$ consists of exactly  $g_{\hat r \hat 4}=1+t_{i,j,k}$ 1-simplices; % - by hypothesis - $K(\e_2,\e_4)=K(r,4)$ consists of exactly one 1-simplex; 
hence, $N(\e_2,\e_4)=N(r,4)$ turns out to be PL-homeomorphic % either  to $\mathbb D^4$ (in the closed case) or to the cone over the boundary $\partial M^4$. 
% of non-empty connected boundary).  
to the boundary connected sum between a handlebody $\mathbb Y^4_{t_{i,j,k}}$  and the cone $v_4 * lkd(v_4, K(\Gamma))$ (which is a cone over $\partial M^4$ in the connected 
boundary case, and a 4-disk in the closed case), where $lkd$ denotes the disjoint link. 
% ; from the singular vertex or to the cone over the boundary $\partial M^4$ (in case of non-empty connected boundary).  

On the other hand, the assumption also implies that both $K(\e_0,\e_3)=K(i,k)$ and $K(\e_1,\e_3)=K(j,k)$ consists of exactly one 1-simplex, while $K(\e_0,\e_1)=K(i,j)$ consists of $g_{\e_2,\e_3,\e_4}=g_{\hat i, \hat  j} \ge 1 $ 1-simplices. 

Moreover, by \cite[Prop. 7.1]{Casali-Cristofori_generalized-genus} (and, in particular, by  \cite[formula (22)]{Casali-Cristofori_generalized-genus}): 

$$ \begin{aligned} 
g_{\widehat{\e_0},\widehat{\e_1},\widehat{\e_3}} \ & =  \ g_{\e_2,\e_4} \ = \ (g_{\widehat{\e_0},\widehat{\e_1}} + g_{\widehat{\e_0},\widehat{\e_3}} -g_{\widehat{\e_0}})  +  \rho_{\e_{\hat0}} \\
\ & = \ (g_{\e_2,\e_3,\e_4} + 1 -1)  +  \rho_{\e_{\hat0}} 
% \ = \ \\ \ & = \ g_{\e_2,\e_3,\e_4}  + \beta_2(M^4)+t_{i,j,k}.% \rho_{\widehat{\e_0}}. 
 \end{aligned} 
 $$  

\par \noindent 
This means that $K(\e_0,\e_1,\e_3)$ contains $\rho_{\e_{\hat0}}$ triangles more than edges in $K(\e_0,\e_1)$. Since 
% $N(\e_0,\e_1,\e_3)$ is PL-homeomorphic either  to $M^4- \mathbb D^4$ (in the closed case) or to $M^4$ (in case of non-empty connected boundary), and since 
each edge in $K(\e_0,\e_1)$ must be face of at least a triangle in  $K(\e_0,\e_1,\e_3)$, then some triangles  in $K(\e_0,\e_1,\e_3)$  may possibly collapse from their ``free" face in $K(\e_0,\e_1)$, giving rise to a 2-dimensional complex $K^\prime(\e_0,\e_1,\e_3)$ consisting of  $\rho_{\e_{\hat0}} +h$ triangles, one edge in $K^\prime (\e_0,\e_3) = K(\e_0,\e_3)$, one edge in $K^\prime(\e_1,\e_3)=K(\e_1,\e_3)$ and $h$ ($1 \le h \le g_{\e_2,\e_3,\e_4}$) edges in $K^\prime(\e_0,\e_1),$
% with   $\beta_2(K^\prime(\e_0,\e_1,\e_3))=\beta_2(M^4)$ and 
such that, for each $i \in \{1, \dots, h\}$, exactly $r_i$ triangles share the same edge in  $K^\prime(\e_0,\e_1)$ (and hence have the same boundary), with $\sum_{i=1}^h r_i = \rho_{\e_{\hat0}} +h.$

Now, it is not difficult to check that, if a ``small'' regular neighbourhood of the $\e_0$-labelled 0-simplex of  $K^\prime(\e_0,\e_1,\e_3)$, together with $h$ (arbitrarily fixed) triangles  of 
$K^\prime(\e_0,\e_1,\e_3)$ having different faces in $K^\prime(\e_0,\e_1),$  is considered as a 0-handle $H^{(0)}=\mathbb D^4,$ then the regular neighbourhoods of the remaining  $\rho_{\e_{\hat0}}$ 
triangles in $K^\prime(\e_0,\e_1,\e_3)$  may be considered - by making use of Lemma \ref{lemma_subgenus} -  as $\rho_{\e_{\hat0}}=\beta_2(M^4)+ t_{i,j,k}$ 2-handles attached on its boundary.

Hence,  
% $N(\e_0,\e_1,\e_3) = H^{(0)} \cup (H^{(2)}_1 \cup \dots \cup H^{(2)}_{\rho_{\widehat{\e_0}}})$, 
 $N(\e_0,\e_1,\e_3) = H^{(0)} \cup (H^{(2)}_1 \cup \dots \cup H^{(2)}_{\beta_2(M^4)+ t_{i,j,k}})$, 
 with $\partial (N(\e_0,\e_1,\e_3)) = \partial (N(\e_2,\e_4))$ PL-homeomorphic to the connected sum between $\partial (\mathbb Y^4_{t_{i,j,k}})= \#_{t_{i,j,k}} (\mathbb S^2 \times \mathbb S^1)$ %(\mathbb S^1 \times \mathbb S^3)
and $\partial M^4$ (in case of non-empty connected boundary). 
%either  to $\mathbb S^3$ (in the closed case) or to $\partial M^4$ (in case of non-empty connected boundary).  

The proof of the general statement is completed by noting that - in virtue of the important result of \cite{Montesinos} and \cite{Laudenbach-Poenaru}, extended in 
\cite{Trace_1982} to the case of simply connected 4-manifolds with connected boundary, as already recalled in Proposition \ref{prop_Trace_3-handles}  -  the boundary identification $\phi$ between $N(\e_0,\e_1,\e_3)$ and $N(\e_2,\e_4)$ is nothing but the attachment of $t_{i,j,k}$  3-handles and, in the closed case, a 4-handle: 
$$ \begin{aligned} 
M^4 \, = \, N(\e_0,\e_1,\e_3) \cup_{\phi} N(\e_2,\e_4) \  = & \ \Big[H^{(0)} \cup \Big(H^{(2)}_1 \cup \dots \cup H^{(2)}_{\beta_2(M^4)+ t_{i,j,k}}\Big)\Big] \cup_{\phi} %\mathbb D^4 
\mathbb Y^4_{t_{i,j,k}}  
\, = \\   = \, & H^{(0)} \cup \Big(H^{(2)}_1 \cup \dots \cup H^{(2)}_{\beta_2(M^4)+ t_{i,j,k}}\Big) \cup  \Big(H^{(3)}_1 \cup \dots \cup H^{(3)}_{t_{i,j,k}}\Big)  \cup     H^{(4)}. \end{aligned} $$

The second part of the statement directly follows form the first one, since $g_{\hat r \hat 4}=1$ means $t_{i,j,k}=0$.
\ \qed 

\medskip

\begin{rem}
{\em
As already pointed out in the introduction, Proposition \ref{handle-decomposition_general} can be equivalently formulated in terms of a relationship between crystallization theory and %(dotted) (undotted) ???? 
framed link representation of compact 4-manifolds. The same relationship, in the inverse direction, is investigated in \cite{Casali-Cristofori_Kirby-diagrams} (and - via a slightly different tool - in  \cite{Casali JKTR2000} and \cite{Casali Compl2004}, too), by describing how to construct a 5-colored graph representing $\bar M^4 (L,c)$ directly from $(L,c)$.   
}\end{rem}   

\bigskip
\smallskip

In virtue of the bijection between colored graphs and colored triangulations described in Section \ref{preliminaries},  
both Theorem  \ref{thm.special-handle-decompositions_closed} and  Theorem  \ref{thm.special-handle-decompositions_boundary} directly follow from Proposition  \ref{handle-decomposition_general}, while the proof of  Theorem  \ref{thm.simple&weaksimple} is an easy consequence: 

%\noindent 
% The proof of Theorem  \ref{framed-link (weak)} is now a direct consequence: 
%  of Proposition  \ref{handle-decomposition}.
\bigskip

\noindent 
{\it Proof of  Theorem   \ref{thm.simple&weaksimple}.} \ 
Let $(\Gamma,\gamma)$ be a weak simple crystallization of the simply-connected compact 4-manifold $M^4$ (with empty or  connected boundary), with respect to the cyclic permutation $\e=(\e_0,\e_1,\e_2,\e_3,\e_4=4)$ of $\Delta_4.$ 
By hypothesis, surely 
$g_{\hat{\e_2} \hat 4}=  g_{\hat{\e_0} \hat{\e_3}} =g_{\hat{\e_1} \hat{\e_3}} = 1$. 
Hence, the second statement of Proposition  \ref{handle-decomposition_general} ensures that the partition $\{\{\e_0,\e_1,\e_3\},\{\e_2,\e_4\}\}$ of $\Delta_4$ induces a handle decomposition of $M^4$ consisting of one 0-handle, $\beta_2(M^4)$ 2-handles and - in the closed case - one 4-handle. This proves that $M^4$ is represented by a (not dotted) framed link  with $\beta_2(M)$ components.
\ \qed

\end{document}